\documentclass[runningheads]{llncs}
\usepackage{graphicx}
\usepackage{latexsym}
\usepackage{amssymb}
\usepackage{amsmath}
\usepackage{mathrsfs}
\usepackage{enumerate}
\usepackage{booktabs}
\usepackage{bm}
\begin{document}
\title{The Logics of Individual Medvedev Frames}
\author{Zhicheng Chen \and Yifeng Ding}
\authorrunning{Z. Chen \and Y. Ding}
\institute{Peking University, Beijing, China}
\maketitle              
\begin{abstract}
Let $n$-Medvedev's logic $\mathbf{ML}_n$ be the intuitionistic logic of 
Medvedev frames based on the non-empty subsets of a set of size $n$, 
which we call $n$-Medvedev frames.
While these are tabular logics, 
after characterizing $n$-Medvedev frames using the property of having at least $n$ maximal points,
we offer a uniform axiomatization of them through a Gabbay-style rule corresponding to this property.
Further properties including compactness, disjunction property, and structural completeness of $\mathbf{ML}_n$
are explored and compared to those of Medvedev's logic $\mathbf{ML}$.
\end{abstract}

\section{Introduction}

In 1962, Y. Medvedev introduced the logic of finite problems \cite{Medvedev1962finiteproblems,Medvedev1966interpretation}, 
also referred to in the literature as Medvedev's logic \textbf{ML}. 
This logic was developed to provide a formal framework 
for Kolmogorov's informal interpretation of intuitionistic logic 
as a calculus of problems \cite{Kolmogorov1932}.

For every positive natural number $n$, 
the $n$-Medvedev frame (up to isomorphism), denoted by $\mathfrak{F}_n$, 
is the Kripke frame 
$$\mathfrak{F}_n=\left\langle \wp^*(n), \supseteq \right\rangle$$
where $\wp^*(X)$ is the set of non-empty subsets of $X$ and 
$n$ is understood set theoretically as $\{0, 1, \dots, n-1\}$. 
The set of Medvedev frames is then $\mathbf{C} = \{\mathfrak{F}_n \mid n>0\}$. 
These frames are essentially finite Boolean algebras 
(viewed as partially ordered structures) without their bottom elements. 
Medvedev’s logic \textbf{ML} is defined as the logic of this class $\mathbf{C}$ 
in the language of intuitionistic propositional logic using Kripke semantics.

By definition \textbf{ML} has the finite model property. 
It is also proved in  \cite{Levin1969syntactic,Maksimova1986disjunction} that 
$\mathbf{ML}$ is a maximal intermediate logic with the disjunction property. 
Moreover, \cite{Prucnal1979twoproblems} reveals that \textbf{ML} is structurally complete. 
As a result, Medvedev’s logic emerges as a notably intriguing and rich intermediate logic. 

In recent work, 
\cite{Holliday2017tolmf} establishes that the class of Medvedev frames is definable by a formula in the extended language of tense logic.
Meanwhile, \cite{Grilletti2022} shows that \textbf{ML} is the logic of finite distributive lattices without the top element,
and \cite{Saveliev2022} shows that \textbf{ML} is the intersection of the logics of finite direct powers of ordinals with the converse ordering taken without the top element.

Despite considerable effort in studying \textbf{ML}, several fundamental logical questions remain unresolved.
For instance, regarding the axiomatizability of \textbf{ML}, 
it is only known that it is not finitely axiomatizable \cite{MaksimovaSS1979}, 
and the question of whether \textbf{ML} is decidable is a long-standing open problem.

In this paper, 
we study the logics $\mathbf{ML}_n$ of each individual $n$-Medvedev frame 
and provide a uniform way to axiomatize them 
with a Gabbay-style rule that intuitively says that 
there are at least $n$ end points in the frame. 
Since these logics are tabular by definition, 
one might wonder if our result is just a corollary of some general theorem about all tabular logics. 
Indeed, it is known that 
all tabular superintuitionstic logics are finitely axiomatizable (Theorem 12.4 of \cite{chagrov1997modal})
using Jankov formulas (Theorem 3.4.27 of \cite{NickBezhanishvilithesis}).
However, to obtain these Jankov formulas for $\mathbf{ML}_n$, 
we have to find all rooted finite Kripke frames of size as large as $2^n$ 
that are not a p-morphic image of a generated subframe of $\mathfrak{F}_{n}$ 
(note that the size of $\mathfrak{F}_n$ is $2^n - 1$). 
While this is a computable task, 
it is unclear whether we can get any syntactic uniformity among these Jankov formulas. 

Our method is as follows. 
First, by improving the characterization of Medvedev frames from \cite{Holliday2017tolmf}, 
we offer a characterization of the $n$-Medvedev frame (among rooted frames) 
using the following three conditions:
\begin{itemize}
    \item (chain${}_{\leq n}$) all strict ascending chains have length at most $n$;
    \item (uni) for any two points $x$ and $y$ seen from $w$ in the frame, 
    there is a third point $z$ seen by $w$ such that 
    the end points seen by $z$ are precisely those seen by either $x$ or $y$;
    \item (end${}_{\geq n}$) there are at least $n$ end points.
\end{itemize}
The conditions (chain${}_{\leq n}$) and (uni) are well known, 
with the first expressible by the standard bounded-depth formula $\boldsymbol{bd_n}$ 
(see for example \cite[p.~278]{vanDalen1986handbook}) 
and the second expressible by the Kreisel-Putnam formula $\boldsymbol{kp}$ \cite{KreiselPutnam1957} 
under the assumption that every point sees at least one end point. 
These two properties do not suffice to characterize $\mathfrak{F}_n$, however, 
and a common supplement is the separability condition (sep) saying that 
two distinct points cannot see exactly the same end points. 
This is, for example, the strategy of \cite{Holliday2017tolmf}. 
The problem of (sep) is that there seems to be no way of capturing it syntactically. 
Our novel observation is that, with (uni) and (chain${}_{\leq n}$), 
we can use (end${}_{\geq n}$) instead (sep) to characterize $\mathfrak{F}_n$. 
While (end${}_{\geq n}$) again cannot be enforced by an axiom, 
it can be by a Gabbay-style rule,
just like irreflexivity can be enforced by the rule 
``to derive $\varphi$, it is enough to derive $(p \land \Box\lnot p) \to \varphi$'' 
\cite{Gabbay1981irreflexivity,Venema1993derivation}. 
Thus, this improved characterization result enables a standard completeness proof through canonical models. 

In addition to the axiomatization result, 
we also discuss some properties of $\mathbf{ML{}_n}$ including compactness, 
disjunction property, and structural completeness, with a comparison to $\mathbf{ML}$. 
In particular, we will show that $\mathbf{ML}$ is not compact.

We introduce notations and the semantic consequence relation 
to be axiomatized in Section \ref{sec:semantic consequence},
and then prove the characterization result in Section \ref{sec:characterization}. 
In Section \ref{sec:axioms} we define our axiomatic system using the new Gabbay-style rule and 
then prove its completeness in Section \ref{sec:completeness}. 
Finally, we discuss the remaining properties in Section \ref{sec:remainings}.

\section{Semantic consequence}
\label{sec:semantic consequence}

Let $PV=\{p_0,p_1,...\}$ be the countably infinite set of propositional variables. 
Let $Form$ collect all the formulas built up from propositional variables and $\perp$, 
using binary connectives $\land,\lor,\to$. 
Negation is defined as usual: $\neg \alpha=(\alpha\to\,\perp)$. 
The order of binding strength of the connectives is $\neg>\lor=\land>\,\to$. 

An intuitionistic Kripke frame is just a non-empty poset, 
while an intuitionistic Kripke model $\left\langle P,\leq, V \right\rangle$ is 
a non-empty poset $\left\langle P,\leq\right\rangle$ endowed with 
an intuitionistic valuation $V: PV \to upset(P, \le)$ 
where $upset(P, \le) = \{X \subseteq P \mid \forall x \in X, y \in P, x \le y \Rightarrow y \in X\}$.
Let $\vDash$ be the forcing relation of the Kripke semantics of intuitionistic propositional logic \textbf{IPL}.

Given a formula $\varphi$, a poset $\mathfrak{P} = (P, \le)$, and a $w \in P$, 
we say that $\mathfrak{P}, w$ validates $\varphi$,
written $\mathfrak{P}, w \vDash \varphi$,
if for any intuitionistic valuation $V$, $\mathfrak{P}, V, w \vDash \varphi$.
By $\Gamma \vDash_{\mathfrak{P}, w} \varphi$, we mean that 
for any intuitionistic valuation $V$, 
if $\mathfrak{P}, V, w \vDash \gamma$ for all $\gamma \in \Gamma$, then $\mathfrak{P}, V, w \vDash \varphi$.
Two more semantic consequence relations can be defined:
\begin{itemize}
    \item $\Gamma \vDash_\mathfrak{P} \varphi$ iff $\Gamma \vDash_{\mathfrak{P}, w} \varphi$ for all $w \in P$;
    \item where $\mathbf{C}$ is a class of non-empty posets, $\Gamma \vDash_\mathbf{C} \varphi$ iff $\Gamma \vDash_\mathfrak{P} \varphi$ for all $\mathfrak{P} \in \mathbf{C}$.
\end{itemize}

\begin{definition}[$n$-Medvedev frame/model]
For every $n > 0$, an \emph{$n$-Medvedev frame} is a poset of the form
$$\mathfrak{F}_k = \langle \wp^*(X), \supseteq \rangle$$
where $X$ is any set of size $n$ and $\wp^*(X)=\wp(X)\setminus\{\emptyset\}$. 
Clearly all $n$-Medvedev frames are isomorphic. For a standard $n$-Medvedev frame, 
let $\mathfrak{F}_k = \langle \wp^*(n), \supseteq \rangle$. 
Here we use the standard set theoretical  definition $n = \{0, 1, \dots, n-1\}$.

An $n$-Medvedev model is an intuitionistic model whose underlying frame is an $n$-Medvedev frame.

Let $\mathbf{C}_{\leq k}=\{\mathfrak{F}_i\mid 1\le i\leq k\}$.
\end{definition}

\begin{definition}[Semantic consequence]
We define the semantic consequence relation of $n$-Medvedev's logic, ${\vDash_n}$, 
by ${\vDash_{\mathfrak{F}_n, n}}$. That is, 
for any $\Gamma \subseteq Form$ and $\alpha \in Form$, 
$\Gamma \vDash_n \alpha\Leftrightarrow\Gamma \vDash_{\mathfrak{F}_n, n} \alpha$.  
\end{definition}

The above definition seems to be only looking at the least element (root) of $\mathfrak{F}_n$, 
but this is only appearance.
\begin{lemma}\label{equivalent definition of vDash_n}
   For any $\Gamma \subseteq Form$ and $\alpha \in Form$, the following are equivalent:
   \begin{enumerate}[(i)]
       \item $\Gamma \vDash_{\mathfrak{F}_n,n} \alpha$;
       \item $\Gamma \vDash_{\mathfrak{F}_n} \alpha$;
       \item $\Gamma \vDash_{\mathbf{C}_{\leq n}} \alpha$.
   \end{enumerate}
\end{lemma}

\begin{proof}
    Clearly $(iii)$ implies $(ii)$ and $(ii)$ implies $(i)$. To show that $(i)$ implies $(iii)$, 
    suppose $\Gamma \not\vDash_{\mathbf{C}_{\le n}} \alpha$.
    Then there are appropriate $\mathfrak{F}_k, V, X$ such that 
    $\mathfrak{F}_k, V, X \vDash \Gamma$ but $\mathfrak{F}_k, V, X \not\vDash \varphi$.
    Now, first take the submodel $(\mathfrak{G}^X, V^X)$ of $(\mathfrak{F}_k, V)$ generated from $X$. 
    It is easy to see that $(\mathfrak{G}^X, V^X)$ is a $|X|$-Medvedev model, and 
    $\mathfrak{G}^X, V^X, X \vDash \Gamma$ while $\mathfrak{G}^X, V^X, X \not\vDash \varphi$.
    In other words, letting $j = |X|$, $\Gamma \not\vDash_{\mathfrak{F}_j, j} \varphi$.
    Now that $n \ge j$, it is well known that there is a surjective p-morphism $f$
    from $\mathfrak{F}_n$ to $\mathfrak{F}_j$ such that $f(n) = j$.
    In fact, for any surjective function $g: n \to j$, the function $f: X \mapsto g[X]$ will do.
    This means $\Gamma \not\vDash_{\mathfrak{F}_n, n} \alpha$.
\end{proof}

\newcommand{\chainn}{(chain${}_{\le n}$)}
\newcommand{\edn}{(end${}_{\ge n}$)}

\section{Characterizing $n$-Medvedev frame}
\label{sec:characterization}

We prove in this section that every rooted poset satisfying (chain${}_{\le n}$), (uni), and (end${}_{\ge n}$) is isomorphic to the $n$-Medvedev frame $\mathfrak{F}_n$.

\begin{definition}
    Let $\mathfrak{P} = \langle P, \leq\rangle$ be a poset. We say that it is \emph{rooted}, when it has a least element, and this element is called its \emph{root}. We call the maximal elements of $\langle P, \leq\rangle$ its \emph{end points}, the set of which is denoted by $End_{\mathfrak{P}}$. For each $w \in P$, let $end_{\mathfrak{P}}(w) = \{u \in End_{\mathfrak{P}} \mid w \leq u\}$, the set of end points above $w$. Subscripts will be dropped when the context is clear.
\end{definition}

\begin{definition}
    Given a poset $\mathfrak{P} = \langle P, \leq\rangle$, we define the following conditions:
    \begin{itemize}
        \item \emph{\chainn} every strictly ascending chain of $\mathfrak{P}$ has at most $n$ elements;
        \item \emph{(uni)} for any $w, u, v \in P$ such that $w \le u, v$, there is $z \in P$ such that $w \le z$ and $end(z) = end(u) \cup end(v)$;
        \item \emph{\edn} $|End| \ge n$.
    \end{itemize}
\end{definition}

\begin{theorem}\label{thm:characterization}
    Let $\mathfrak{P} = (P, \le)$ be a rooted poset satisfying \chainn, (uni), and \edn. Then $\mathfrak{P}$ is isomorphic to $\mathfrak{F}_n$.
\end{theorem}
\begin{proof}
    Let $r$ be the root of $\mathfrak{P}$. First, note that for any $w \in P$, $end(w) \not= \varnothing$, as otherwise there will be an infinite ascending chain starting from $w$. Now we claim $|End| \le n$. Suppose not. Then take $n+1$ distinct elements $e_1, e_2, \dots, e_{n+1}$ from $End$. Now let $u_1 = e_1$ and inductively let $u_{i+1}$ be the element promised by (uni) when applied to $r$, $u_i$, and $e_{i+1}$. Then we see that for any $i$ from $1$ to $n+1$, $end(u_i) = \{e_1, e_2, \dots, e_i\}$. This means for any $i = 1, \dots, n$, $u_i < u_{i+1}$ since $u_i \le u_{i+1}$ by construction from (uni) and they cannot be equal because $end(u_i) \subsetneq end(u_{i+1})$. Then we obtain a strictly ascending chain of length $n+1$, contradicting \chainn. Together with \edn, we have shown that $|End| = n$. 
    Clearly, all that is left to be shown is that 
    $end$ as a function from $P$ to $\wp^*(End)$ is a bijection and 
    for any $u, v \in P$, $u \le v$ iff $end(u) \supseteq end(v)$.

    If $u \le v$, then trivially $end(u) \supseteq end(v)$.
    Conversely, suppose $u \not\le v$ yet $end(u) \supseteq end(v)$.
    Then apply (uni) to $r$, $u$, and $v$ to obtain a $z$ such that 
    $z \le u, v$ and $end(z) = end(u) \cup end(v) = end(u)$. 
    Since $u \not\le v$, $z < u$.
    Now we obtain a strictly ascending chain of length $n+1$ since 
    there must be a strictly ascending chain of length $|end(u)|$ starting from $u$, and 
    there must be a strictly descending chain of length $1 + (n - |end(u)|)$ starting from $z$.
    More specifically, let $k = |end(u)|$ and 
    list $end(u) = \{e_1, \dots, e_k\}$ and $End \setminus end(u) = \{e_{k+1}, \dots, e_n\}$. 
    Define $a_1 = e_1$ and inductively pick $a_{i+1}$ by applying (uni) to $u$, $a_i$, and $e_{i+1}$.
    Similarly, define $b_k = z$ and inductively pick $b_{i+1}$ by applying (uni) to $r$ (the root), $b_i$, and $e_{i+1}$.
    Then we can observe that $b_n < b_{n-1} < \dots < b_{k+1} < z < u < a_{k-1} < \dots < a_1$ is a strictly ascending chain of length $n+1$, contradicting \chainn.
    This also shows that $end: P \to \wp(End)$ is injective, 
    since if $end(u) = end(v)$ then $u \le v$ and $v \le u$. 

    That $end: P \to \wp(End)$ is surjective is an easy consequence of 
    iteratively applying (uni) to $r$ and any subset of $End$. 
\end{proof}

\newcommand{\redn}{\ensuremath{\boldsymbol{Ed}_n}}

\section{Axiomatization and soundness}
\label{sec:axioms}
As we mentioned in the introduction, 
the key to our uniform axiomatizations of the logics of individual Medvedev frames 
is a Gabbay-style rule that intuitively says there are at least $n$ end points. 

\begin{definition}
    Let \redn be the following rule:
    \begin{equation}
    \begin{aligned}
        &\alpha \to (\beta \lor \bigvee_{i = 1}^n \lnot (p_i \land \bigwedge_{j \not= i} \lnot p_j)) \\ 
        \midrule     
        &\alpha \to \beta
    \end{aligned}
    \tag{\redn}
    \end{equation}
    where $p_1, \dots, p_n$ must be distinct propositional variables not occurring in $\alpha$ or $\beta$. 
    For convenience, we write $\lambda_i$ for $p_i \land \bigwedge_{j \not= i} \lnot p_j$.
\end{definition}

\begin{lemma}\label{frame correspondence}
Let $\mathfrak{P}=\langle P,\leq\rangle$ be a poset with root $r$. We say $u, v \in P$ are \emph{incompatible} if there is no $z$ above both $u$ and $v$. Then, the following are equivalent: 
\begin{enumerate}
    \item[(1)] $\mathfrak{P}, r$ validates \redn. That is, whenever $\mathfrak{P}, r$ validates $
        \alpha \to (\beta \lor \bigvee_{i = 1}^n \lnot \lambda_i)$ where $p_1, \dots, p_n$ are distinct propositional variables not occurring in $\alpha$, $\mathfrak{P}, r$ also validates $\alpha \to \beta$.
    \item[(2)] $\mathfrak{P}$ has at least $n$ pairwise incompatible elements.
\end{enumerate}
If additionally for any $w \in P$, $end(w) \not= \varnothing$, then the above is equivalent to that $|End| \ge n$. 
\end{lemma}
\begin{proof}
    (1) $\Rightarrow$ (2): Suppose $\mathfrak{F}, r$ validates \redn. 
    Consider the case where $\alpha = \top$ and $\beta = \bot$. 
    Since $\mathfrak{F}, r$ does not validate $\top \to \bot$, 
    $\mathfrak{F}, r$ must not validate $\top \to (\bot \lor \bigvee_{i = 1}^n\lnot\lambda_i)$, 
    which is just $\bigvee_{i = 1}^n\lnot\lambda_i$.
    Thus there is an intuitionistic valuation $V$ such that for any $i = 1, \dots, n$, 
    $\mathfrak{F}, V, r \not\vDash \lnot\lambda_i$, meaning that 
    for some $u_i \in P$, $\mathfrak{F}, V, u_i \vDash \lambda_i$.
    Since $\lambda_i$ is $p_i \land \bigwedge_{j \not= i} \lnot p_j$, 
    $u_1, \dots, u_n$ are pairwise incompatible.
    
    (2) $\Rightarrow$ (1): 
    Suppose $\mathfrak{P}$ has $n$ pairwise incompatible elements $u_1, \dots, u_n$.
    To show that $\mathfrak{P}, r$ validates \redn, 
    suppose $\mathfrak{P}, r$ does not validate $\alpha \to \beta$.
    That is, we have an intuitionistic valuation $V$ such that $\mathfrak{P}, V, r \vDash \alpha$,
    but $\mathfrak{P}, V, r \not\vDash \beta$.
    Let $V'$ be the valuation that are identical to $V$ except that 
    $V'(p_i) = {\uparrow\!u_i} := \{z \in P \mid z \ge u_i\}$. 
    Since $u_i$s are pairwise distinct, $\mathfrak{F}, V', u_i \vDash \lambda_i$. 
    Thus, $\mathfrak{P}, V', r \not\vDash \lnot\lambda_i$ for each $i = 1, \dots, n$.
    Since $V$ and $V'$ are identical except for $p_1, \dots, p_n$ that do not occur in $\alpha$ or $\beta$, 
    $\mathfrak{F}, V', r$ continues to force $\alpha$ but not $\beta$. 
    So $\mathfrak{F}, r$ invalidates $\alpha \to (\beta \lor \bigvee_{i = 1}^n \lnot\lambda_i)$ using $V'$.

    Finally, if for any $w \in P$, $end(w) \not= \varnothing$, 
    then clearly there are $n$ pairwise incompatible elements iff $|End| \ge n$.
\end{proof}

Now we define our axiomatization.

\begin{definition}
    Let $\mathbf{ML}_n$ be the smallest superintuitionistic logic that contains all instances of
    \begin{align*}
        &\boldsymbol{kp}\ \qquad \quad (\neg p \rightarrow q \vee r) \to (\neg p \rightarrow q) \vee(\neg p \rightarrow r)\\
        &\boldsymbol{bd_n}\qquad\quad p_n \vee(p_n \rightarrow(p_{n-1} \vee(p_{n-1} \rightarrow(...(p_1 \vee(p_1 \rightarrow \perp
        ))...))))
    \end{align*}
    and is further closed under \redn. 
    By $\Gamma \vdash_{n} \varphi$, 
    we mean that there are finitely many formulas $\gamma_1, \dots, \gamma_n \in \Gamma$ such that 
    $\bigwedge_i \gamma_i \to \varphi$ is provable in $\mathbf{ML}_n$.
\end{definition}

Since $\mathbf{ML}_n$ is a superintuitionistic logic, 
$\vdash_n$ enjoys the basic properties of the intuitionistic consequence relation.
But because of \redn, it also has the following property:
\begin{lemma}\label{lem:rednconsequence}
    If $p_1, \dots, p_n$ are distinct propositional variables not occurring in $\Gamma \cup \{\varphi\}$ and $\Gamma \vdash_n \varphi \lor \bigvee_{i = 1}^n \lnot\lambda_i$, $\Gamma \vdash_n \varphi$.
\end{lemma}

We end this section with the soundness of $\mathbf{ML}_n$.
\begin{theorem}
    For any $\varphi$ provable in $\mathbf{ML}_n$, $\varphi$ is validated by $\mathfrak{F}_n, n$.
\end{theorem}
\begin{proof}
    The validity of $\boldsymbol{kp}$ is well-known.
    It is also well known that $\boldsymbol{bd}_n$ is valid at a point in a poset iff 
    all strictly ascending chain starting from that point has length at most $n$, 
    which is true at the root of $\mathfrak{F}$. 
    Finally, since $\mathfrak{F}_n$ has exactly $n$ end points, by Lemma \ref{frame correspondence},
    \redn\ preserves validity at $\mathfrak{F}_n, n$.
\end{proof}

\section{Completeness}
\label{sec:completeness}
We first define the canonical model for $\mathbf{ML}_n$.
\begin{definition}
For any $\Gamma \subseteq Form$, 
\begin{enumerate}
  \item it is \emph{closed}, if for each  $\varphi\in Form$, $\Gamma\vdash_n\varphi$ implies $\varphi\in\Gamma$;
  \item it is \emph{consistent}, if $\Gamma\nvdash_n\perp$.
  \item it is \emph{prime}, if for any $\varphi,\psi\in Form$, $\varphi\vee\psi\in\Gamma$ implies either $\varphi\in\Gamma$ or $\psi\in\Gamma$.
\end{enumerate}
Then, define 
\begin{itemize}
  \item $P^c=\{\Gamma \subseteq Form \mid \Gamma \text{ is consistent, closed and prime}\}$;
  \item $\mathfrak{P}^c = (W^c, \subseteq)$;
  \item $V^c : PV \to \wp(W^c)$ such that for each $p \in PV$, $V^c(p) = \{\Gamma\in W^c \mid p\in\Gamma\}$.
\end{itemize}
Then define the canonical model $\mathcal{M}^c=\langle \mathfrak{P}^c, V^c \rangle$.
\end{definition}

Since $\vdash_n$ is a superintuitionistic consequence relation, the following lemmas are standard.
\begin{lemma}[Lindenbaum’s lemma]\label{Lindenbaum}
For any $\Gamma \subseteq Form$ and $\alpha \in Form$, if $\Gamma \nvdash_n \alpha$, then there exists $\Phi \in W^c$, such that $\Gamma \subseteq \Phi$ and $\alpha \notin \Phi$.
\end{lemma}

\begin{lemma}[Existence lemma]\label{existence lemma}
    For any $\Gamma \subseteq Form$ and $\alpha, \beta \in Form$ , 
    if $\Gamma \nvdash_n \alpha \rightarrow \beta$,
    then there exists $\Phi \in P^c$ such that 
    $\Gamma \subseteq \Phi$ and $\alpha \in \Phi$ and $\beta \notin \Phi$.
    In particular, if $\Gamma \not\vdash_n \lnot\alpha$, 
    then there is $\Phi \in P^c$ such that $\Gamma \cup \{\alpha\} \subseteq \Phi$.
\end{lemma}

\begin{lemma}[Truth lemma]
For any $\varphi\in Form$ and $\Gamma\in W^c$,
\[\mathcal{M}^c,\Gamma\vDash\varphi \Leftrightarrow \varphi\in\Gamma.\]
\end{lemma}

It is also well known that $\boldsymbol{bd}_n$ and $\boldsymbol{kp}$ are canonical
(see e.g. Theorem 5.16 of \cite{chagrov1997modal}).
Thus, we know the following of $\mathfrak{P}^c$.
\begin{lemma}
    $\mathfrak{P}^c$ satisfies \chainn. Thus, for each $\Phi\in P^c$, $end(\Phi) \not= \varnothing$.
\end{lemma}
Note that even without $\boldsymbol{bd}_n$, $end(\Phi)$ is still non-empty for any $\Phi \in P^c$,
since we can always extend $\Phi$ to a maximally consistent set.
\begin{lemma}
    $\mathfrak{P}^c$ satisfies (uni).
\end{lemma}
\begin{proof}
    While the validity of $\boldsymbol{kp}$ only corresponds to the following property on arbitrary non-empty poset $(P, \le)$:
    \begin{align*}
    \forall x,y,z &(x \le y \land x \le z \land y \not\le z \land z \not\le y \to\\
    &\exists u \in {\uparrow\! x} (u \le y\land u \le z \land (\forall v \in {\uparrow\! u} (\exists w \in {\uparrow\! v}((y \le w \lor z \le w)))))),
    \end{align*}
    (uni) follows from this if $end(w) \not= \varnothing$ for every $w \in P$. 
    But as we commented above, $\mathfrak{P}^c$ is indeed such that for every $\Phi \in P^c$, 
    $end(\Phi) \not= \varnothing$.
\end{proof}

To make use of \redn, we must pay special attention to avoid using some propositional variables.
\begin{lemma}\label{canonical end_n}
    If $\Gamma \nvdash_n \alpha$ and 
    $p_1, \ldots, p_n$ are distinct propositional variables not occurring in $\Gamma\cup\{\alpha\}$, 
    then there exists $\Phi \in P^c$ such that 
    $\Gamma \subseteq \Phi$, $\alpha \notin \Phi$, and $|end(\Phi)|\geq n$.
\end{lemma}

\begin{proof}
    Using \redn\ and Lemma \ref{lem:rednconsequence}, 
    $\Gamma \nvdash_n \alpha \vee \neg\lambda_1\vee...\vee\neg\lambda_n$.
    Then by lemma \ref{Lindenbaum}, there exists $\Phi \in P^c$ s.t. 
    $\Gamma \subseteq \Phi$ and $\alpha\vee\neg \lambda_1\vee...\vee\neg\lambda_n \notin \Phi$.
    Since $\Phi$ is closed, none of the disjuncts is in $\Phi$.
    By Lemma \ref{existence lemma} and using maximally consistent sets,
    for each $1\leq i\leq n$, there is a $\Theta_i\in end(\Phi)$ s.t. $\lambda_i\in\Theta_i$.
    Since $\lambda_i$s are pairwise inconsistent, $\Theta_i$s are pairwise distinct.
    Thus $|end(\Phi)| \ge n$.
\end{proof}
Now we can prove a conditioned completeness theorem.
\begin{lemma}\label{lem:conditional completenss}
    If there are distinct $p_1, \ldots, p_n\in PV$ not occurring in $\Gamma\cup\{\alpha\}$, then $\Gamma \nvdash_n \alpha$ implies $\Gamma \nvDash_n \alpha$.
\end{lemma}
\begin{proof}
    By lemma \ref{canonical end_n}, there is $\Phi \in W^c$ such that
    $\Gamma \subseteq \Phi$, $\alpha \notin \Phi$, and $|end(\Phi)|\geq n$. 
    Moreover, by the Truth Lemma, $\mathcal{M}^c,\Phi\vDash\Gamma$ and $\mathcal{M}^c,\Phi\nvDash\alpha$.
    Let $\mathcal{M}^\Phi = \langle \mathfrak{P}^\Phi, V^\Phi \rangle$ be 
    the submodel of $\mathcal{M}^c$ generated from $\Phi$.
    Then $\mathcal{M}^\Phi,\Phi \vDash \Gamma$ and $\mathcal{M}^\Phi,\Phi \nvDash \alpha$.

    Since $\mathfrak{P}^c$ satisfies \chainn and (uni), $\mathfrak{P}^\Phi$ does so as well. 
    Since $|end(\Phi)| \geq n$, $\mathfrak{P}^\Phi$ satisfies \edn. 
    By Theorem \ref{thm:characterization}, $\mathfrak{P}^\Phi$ is isomorphic to $\mathfrak{F}_n$.
    Thus $\Gamma \nvDash_n \alpha$.
\end{proof}

\begin{theorem}[Strong completeness theorem]
    For any $\Gamma \cup \{\alpha\} \subseteq Form$,
    \[\Gamma\,\vDash_n\alpha \Rightarrow\Gamma\,\vdash_n\alpha.\]
\end{theorem}
\begin{proof}
    The difficulty lies in the case where we cannot find $n$ distinct propositional variables
    not occurring in $\Gamma\cup\{\alpha\}$.

    Since we have stipulated that $PV$ is countably infinite, fix a injection $\iota: PV \to PV$
    such that there are $n$ distinct $p_1, \dots, p_n \in PV$ not in the range of $\iota$.
    Treat $\iota$ as a substitution and extend it to a function from $Form$ to $Form$. 
    Then it is easy to observe that $\Gamma \vDash_n \alpha$ iff 
    $\iota[\Gamma] \vDash_n \iota(\alpha)$ (by shuffling the valuations). 
    For $\iota[\Gamma]$ and $\iota(\alpha)$, we use Lemma \ref{lem:conditional completenss} and
    obtain a proof $\pi$ of $\bigwedge_{i = 1}^k \iota(\gamma_i) \to \iota(\alpha)$
    where $\gamma_1, \dots, \gamma_k \in \Gamma$.
    Now, the proof $\pi$ is finite, so let $Q$ collect the propositional variables used in $\pi$.
    Then there is also an injection $\iota': PV \to PV$ such that 
    for any $q \in Q$, 
    if $q$ is in the range of $\iota$,  then $\iota'(q) = \iota^{-1}(q)$, and 
    otherwise $\iota'(q)$ is some new propositional variable not in $\iota^{-1}[Q]$. 
    This $\iota'$ is like $\iota^{-1}$ when applied as a substitution to formulas in $\pi$ 
    and also the range of $\iota$.
    Thus, when applying $\iota'$ to formulas in $\pi$, we obtain a proof of $\bigwedge_{i = 1}^k\gamma_i \to \alpha$,
    meaning that $\Gamma \vdash_n \alpha$.
\end{proof}

\section{Other properties and comparison with Medvedev's logic}
\label{sec:remainings}

In the remainder of this paper, we look at some properties of $n$-Medvedev's logic with comparison to Medvedev's logic.
Let $\mathbf{ML}$ be Medvedev's logic.

We first observe that 
\begin{proposition}
    $\mathbf{ML} \subsetneq \dots \subsetneq \mathbf{ML}_{n+1} \subsetneq \mathbf{ML}_n \subsetneq \mathbf{ML}_{n-1} \subsetneq \dots \subsetneq \mathbf{ML}_1$ and $\mathbf{ML} = \bigcap_n \mathbf{ML}_n$.
\end{proposition}
By the well known result of Levin \cite{Levin1969syntactic} and Maksimova \cite{Maksimova1986disjunction}, \textbf{ML} is a maximal logic with the disjunction property. So, it follows from $ML\subset ML_n$ that:
\begin{proposition}
    $\textbf{ML}_n$ does not enjoy the disjunction property.
\end{proposition}

The strong completeness of $n$-Medvedev's logic implies its compactness, 
though this is also follows from first-order definability of the class of $n$-Medvedev frames.

However, it is not the case for Medvedev's logic:
\begin{proposition}
    Let $\mathbf{C} = \{\mathfrak{F}_n \mid n \ge 1\}$. 
    Then the semantic consequence $\vDash_\mathbf{ML}$ of Medvedev's logic 
    defined by $\vDash_\mathbf{C}$ is not compact, 
    as illustrated by the following counterexample, where $p_0$ is not used in any $\boldsymbol{bd}_i$:
    $$\{\boldsymbol{bd}_i\to p_0 \mid i > 0\}\vDash_{\textbf{ML}} p_0,$$ 
    but $\Gamma'\nvDash_{\textbf{ML}} p_0$ for any finite 
    $\Gamma'\subseteq\{\boldsymbol{bd}_i\to p_0\mid i > 0\}$. 
\end{proposition}
\begin{proof}
    On the one hand, 
    let $\mathcal{M},r$ be an arbitrary pointed Medvedev model with 
    $\mathcal{M},r\vDash\{\boldsymbol{bd}_i \to p_0 \mid i > 0\}$. 
    Let $k$ be the number of end points of $\mathcal{M}$. 
    Then the underlying frame of $\mathcal{M}$ satisfies (Chain${}_{\leq k}$), 
    and thus $\mathcal{M},r\vDash\boldsymbol{bd}_k$. 
    Together with $\mathcal{M},r\vDash\boldsymbol{bd}_k \to p_0$, 
    it follows that $\mathcal{M},r\vDash p_0$.

    On the other hand, 
    let $\Gamma'\subseteq\{\boldsymbol{bd_i}\to p_0\mid i\in\omega^*\}$ be finite. 
    Clearly $\,\nvDash_{\textbf{ML}} p_0$. 
    So, we may assume that $\Gamma'$ is not empty. 
    Note that $\boldsymbol{bd_i}\vDash_{\textbf{ML}}\boldsymbol{bd_j}$ for any $i < j$, 
    which means $\boldsymbol{bd}_i \to p_0 \vDash_{\mathbf{ML}} \boldsymbol{bd}_j \to p_0$ for any $i > j$. 
    So we only need to prove $\boldsymbol{bd}_{i} \to p_0\nvDash_{\textbf{ML}} p_0$ with 
    the largest $i$ so that $\boldsymbol{bd}_{i} \to p_0 \in \Gamma'$.
    But a countermodel can be defined easily on the $(i+1)$-Medvedev frame.
\end{proof}

Prucnal \cite{Prucnal1979twoproblems} proved that Medvedev's logic is structurally complete,
in the sense that for any $\varphi,\psi\in Form$, 
$\varphi\to\psi\in \mathbf{ML}$ if for all substitutions $\sigma$, 
$\sigma(\varphi)\in \mathbf{ML}$ implies $\sigma(\psi)\in \mathbf{ML}$.
By adapting the proof in \cite{Prenosil2024levinsprucnalstheoremsmedvedevs} 
of the structural completeness of Medvedev's logic, 
we also show the structural completeness of $n$-Medvedev's logic. 

\begin{proposition}
    $n$-Medvedev's logic is structurally complete.
\end{proposition}

\begin{proof}
    We prove the contrapositive. 
    Suppose $\varphi\to\psi\notin \mathbf{ML}_n$. 
    Then $\varphi\nvDash_{n}\psi$. 
    By definition, there is an intuitionistic valuation $V$ on $\mathfrak{F}_n$ such that 
    $\mathfrak{F}_n,V,n \vDash \varphi$ and $\mathfrak{F}_n,V,n\nvDash \psi$.
    Let $\alpha_{\{0\}}$,...,$\alpha_{\{n-1\}}\in Form$ and 
    $V_0$ be an intuitionistic valuation on $\mathfrak{F}_n$ such that 
    for any $0\leq i,j < n$,
    \begin{enumerate}[(i)]
        \item $\vdash_{\textbf{IPL}}\neg(\alpha_{\{i\}}\land \alpha_{\{j\}})$ if $i\neq j$;
        \item $\vdash_{\textbf{IPL}}\neg\neg(\alpha_{\{0\}}\lor...\lor \alpha_{\{n-1\}})$;
        \item $\mathfrak{F}_n,V_0,\{j\} \vDash \alpha_{\{i\}}\Leftrightarrow i=j$.
    \end{enumerate}
    (For an example of $\alpha_{\{0\}}$,...,$\alpha_{\{n-1\}}$ and $V_0$, see \cite{Prenosil2024levinsprucnalstheoremsmedvedevs}.)
    
    For any $I\in\wp^*(n)$, we define $\alpha_I = \neg\neg\bigvee_{i\in I}\alpha_{\{i\}}$. 
    Then it is easy to see that for any $J\in\wp^*(n)$, 
    \[\mathfrak{F}_n,V_0,J \vDash \alpha_{I}\Leftrightarrow I\supseteq J.\]
    In a sense, formula $\alpha_I$ characterizes the principle upset ${\uparrow\!I}$ in 
    $\mathfrak{F}_n=(\wp^*(n),\supseteq)$. 
    Since $\wp^*(n)$ is finite, 
    every upset $S$ in $upset(\mathfrak{F}_n)$ can be characterized by \[\alpha_S = \bigvee_{I\in S}\alpha_I .\] 

    Now, define the substitution $\sigma$ : for each $p\in PV$, $\sigma(p)=\alpha_{V(p)}$. 
    Then it not hard to show, by induction, that for any $\beta\in Form$, for any $J\in\wp^*(n)$, 
    \[\mathfrak{F}_n,V_0,J \vDash \sigma(\beta)\Leftrightarrow \mathfrak{F}_n,V,J \vDash\beta.\]

    Since $n\in\wp^*(n)$ and $\mathfrak{F}_n,V,n\nvDash \psi$, we have $\mathfrak{F}_n,V_0,n\nvDash \sigma(\psi)$. So $\nvDash_n \sigma(\psi)$.
    
    It remains to prove that $\vDash_n \sigma(\varphi)$. 
    By $\mathfrak{F}_n,V,n\vDash \varphi$, we get $\mathfrak{F}_n,V_0,n\vDash \sigma(\varphi)$.
    Due to the special form of formulas produced by substitution $\sigma$, 
    as well as the nice properties (i),(ii) of $\alpha_{\{0\}}$,...,$\alpha_{\{n-1\}}$, 
    it turns out that from $\mathfrak{F}_n,V_0,n\vDash \sigma(\varphi)$, 
    one can in fact prove $\sigma(\varphi)\in \mathbf{ML}$. 
    (We refer the reader to \cite{Prenosil2024levinsprucnalstheoremsmedvedevs} for details.)
    Since $\mathbf{ML} \subseteq \mathbf{ML}_n$, we have $\vDash_n \sigma(\varphi)$.\qed
\end{proof}

%
%
%
\bibliographystyle{splncs04}
\bibliography{main}
\end{document}